\documentclass[12pt,english]{article}
\usepackage{cite,amsfonts,amssymb,graphics} 

\begin{document} 

\title{\bf Numerical solutions to integrodifferential equations
which interpolate heat and wave equations
\thanks{~Extended version of the
talk given by P.R. at Second International Conference 
of Applied Mathematics, Plovdiv, Bulgaria, August 12-18, 2005} }

\author{\sf Piotr Rozmej$\;^1$ and Anna Karczewska$\;^2$  \\[2mm]
$^1$ Institute of Physics, University of Zielona G\'ora\\
  ul. Szafrana 4a, 65-246 Zielona G\'ora, POLAND\\
  e-mail:  P.Rozmej@if.uz.zgora.pl\\[2mm]
$^2$ Department of Mathematics,
 University of Zielona G\'ora\\
 ul. Szafrana 4a, 65-246 Zielona G\'ora, Poland\\
 e-mail: A.Karczewska@im.uz.zgora.pl\\
}

\maketitle

\def\thefootnote{}

\begin{abstract}
In the paper we study some numerical solutions to
Volterra equations which interpolate heat and wave equations.
We present a scheme for construction of approximate numerical solutions
for one and two spatial dimensions. Some solutions to the stochastic 
version of such equations (for one spatial dimension)
are presented as well.
\end{abstract}
\vspace{3mm}

\vfill\newpage

\section{Introduction}

We consider the following integrodifferential equation (Volterra type)
\vspace{1mm}
\begin{equation} \label{e1}
f(x,t) = g(x) + \int_0^t a(t-s)\,A\,f(x,s)\, ds \;,
\end{equation}
where $A$ is Laplacian and 
$\displaystyle a(t)\!=\! \frac{t^{\alpha-1}}{\Gamma(\alpha)},~\Gamma$
is the gamma function,
$1 \leq \alpha \leq 2$, $x\in\mathbb{R}^d, ~t\geq 0.$ 
The equation
(\ref{e1}) was considered in context of 
 the heat con\-duction with memory \cite{GuPi,Mill}.

For particular cases $\alpha=1$ and 
$\alpha=2$ the equation (\ref{e1}), 
after taking the first and the second time derivative, 
becomes the heat and the wave equation, 
respectively. 
For $1< \alpha < 2$ the equation (\ref{e1}) 
interpolates the heat and the wave equations.
The equation (\ref{e1}) was discussed extensively by Fujita \cite{Fuji}
and Schneider \& Wyss \cite{SchW}. Fujita \cite{Fuji}
has found the analytical form 
of solutions $f(x,t)$ to  (\ref{e1}) in terms of 
{\em resolvents} or {\em fundamental solutions}  $S(t)$.

A stochastic version of the equation (\ref{e1})
\begin{equation}\label{e4}
f(x,t) = g(x) + \int_0^t a(t-s)\,A\,f(x,s)\, ds +W(x,t)\,,
\end{equation}
where ~$W$  is some stochastic process has been studied in 
\cite{KaIN} and \cite{KaRo}.

Within the resolvent approach the {\it mild solution} to 
(\ref{e4}) is given in the form:
\begin{equation}\label{e5}
f(x,t) = S(t)f(x,0) + \int_0^t S(t-\tau)dW(x,\tau)\,,
\end{equation} 
where the operator  $S(t)$ is the resolvent (fundamental solution)
to the equation (\ref{e1}), i.e.\  $f(x,t) = S(t)f(x,0)$.
The resolvent $S(t)$ found by Fujita \cite{Fuji} is given by the 
formula
\begin{equation}\label{e5a}
 (S(t)f)(x) = \int_{-\infty}^{\infty} \phi_\alpha(t,x-y)f(y)dy,
 \quad t\geq 0,\quad x\in \mathbb{R}\,,
\end{equation}
where 
\begin{eqnarray}\label{e5b}
\phi_\alpha(t,x) &=&e^{-x^2/4t}/\sqrt{4\pi t} 
\quad\mathrm{for}\quad \alpha=1 \quad\quad \mathrm{and}
\nonumber \\
 \phi_\alpha(t,x) &=& 
\frac{1}{2}(\delta_0(t-x)+\delta_0(t+x))
\quad\mathrm{for}\quad \alpha=2, 
\end{eqnarray}
($\delta_0(x)$--Dirac's $\delta$-function).
For $1< \alpha <2$, the analytical form of $\phi_\alpha(t,x)$ 
is given in terms of inverse Fourier transform 
of Mittag-Leffler function
$ML_{\alpha}(z)$ \cite{Fuji,Pr} and a direct calculation of both
solutions to (\ref{e1}) and resolvents becomes very difficult.
It seems that obtaining approximate numerical solutions may be more 
practical. 
 
The aim of the paper is to construct: 
\begin{itemize}
\item approximate numerical solutions to equations (\ref{e1}), 
 (\ref{e4}) (deterministic and stochastic) for $d=1$,
\item  numerical solutions to equation (\ref{e1}) for $d=2$.
\end{itemize}
The existing analytical solutions to (\ref{e1}) for $d=1$
will serve as a reference to control the quality of the 
numerical approximation.

For arbitrary $1\leq \alpha \leq 2$, the resolvent operator $S(t)$
for (\ref{e1}) does not possess a semigroup property. 
Hence, the time evolution from 0 to $t$ can not be divided into 
smaller steps and has to be calculated in one step.
Therefore, the Galerkin method for numerical approach is a reasonable
choice.

The paper is organized as follows. In section 2 the Galerkin method 
for solving (\ref{e1}) with one spatial dimension is presented. 
The numerical solutions for $\alpha=1$ and $\alpha=2$ are 
compared to existing analytical ones. Examples of numerical solutions 
for stochastic cases with a simple stochastic process are presented as
well.
In section 3 the Galerkin method for two spatial dimensions is
presented. Several results of numerical solutions to (\ref{e1})
for different $\alpha$ are shown, too.

\section{Galerkin method, case $d=1$}

In Galerkin method one introduces a complete set of orthonormal 
functions $ \{\phi_j\}, ~j=1,2,\ldots ,\infty $ on the interval $[0,t]$,
spanning a Hilbert space $H$. 
Then the approximate solution is postulated as an expansion of the
unknown true solution in the subspace $H_n$ spanned by $n$ first basis 
functions $ \{\phi_k\}, ~k=1,2,\ldots ,n $
\begin{equation}  \label{e6a}
f_n(x,t)=\sum_{k=1}^n c_k(x) \phi_k(t)\;.
\end{equation} 
Inserting (\ref{e6a}) into (\ref{e1})  we obtain
\begin{equation} \label{e6b}
f_n(x,t) = f(x,0) + \int_0^t  a(t-s)\frac{d^2}{dx^2}
\,f_n(x,s) ds +\varepsilon_n(x,t) \;,
\end{equation}
where the function $\varepsilon_n(x,t)$ represents 
the approximation error. 
From (\ref{e6a}) and (\ref{e6b}) we have
\begin{eqnarray} \label{e6c}
\varepsilon_n(x,t) & = & f_n(x,t) - f(x,0) - \int_0^t  
a(t-s)\frac{d^2}{dx^2}\,f_n(x,s) ds \\ & = &
\sum_{k=1}^n c_k(x)\phi_k(t) - f(x,0) -\int_0^t 
a(t-s)\frac{d^2}{dx^2}\sum_{k=1}^n c_k(x)\phi_k(s) \,ds
\;. \nonumber
\end{eqnarray}
Requirement that the error function $\varepsilon_n(x,t)$
 has to be orthogonal to the subspace $H_n$, 
 $(\phi_j(t), \varepsilon_n(x,t))=0$, for $j=1,2,\ldots ,n$,
 leads to the set 
 of coupled differential equations for the coefficient
 functions $c_j(x)$ 
 \begin{equation} \label{e7}
 g_j(x) = c_j(x) - \sum_{k=1}^n  a_{jk}\frac{d^2 c_k(x)}{dx^2}  \;,
\end{equation} 
where
\begin{equation} \label{e8}
 a_{jk}  = \int_0^t\! 
\phi_j(\tau)\left[ 
\int_0^\tau \! a(\tau-s)\phi_k(s)ds \right] d\tau, ~~
\mbox{(in general}~~ a_{jk} \neq a_{kj})
\end{equation} 
and 
\begin{equation} \label{e8a}
 g_j(x) = \int_0^t f(x,0)\phi_j(\tau)d\tau = f(x,0)\int_0^t 
 \phi_j(\tau)d\tau \,.
\end{equation}

Discretizing second derivative (Laplacian) one obtains (\ref{e7}) as: 
\begin{equation} \label{e9}
g_j(x_i) = c_j(x_i) + \frac{1}{h^2} \sum_{k=1}^n a_{jk}\, 
 [-c_k(x_{i-1})+2 c_k(x_i) - c_k(x_{i+1})]
\end{equation} 
with ~$h\!=\! x_{i}-x_{i-1}$ and $j=1,2,\ldots,n,~i=1,2,\ldots,m$.

The set (\ref{e9}) can be written in matrix form: 
$\underline{\mathcal{A}}\, \underline{c} = \underline{g}$, where  
$\underline{c}$ and $\underline{g}$ are 
$(N=n\cdot m)$-dimensional vectors and matrix 
$\underline{\mathcal{A}}$ has a block form
\begin{equation} \label{e10}
\underline{c} = \left( \begin{array}{c} C_1 \\ C_2 \\ \vdots \\ C_n
  \end{array}\right),
  \hspace{2ex}
  \underline{g} = \left( \begin{array}{c} G_1 \\ G_2 \\ \vdots \\ G_n
  \end{array}\right),
  \hspace{2ex}
  \underline{\mathcal{A}} =  \left( \begin{array}{ccc}
 \left[A_{11}\right] & \ldots & \left[A_{1n}\right] \\ 
 \left[A_{21}\right] & \ldots & \left[A_{2n}\right] \\ 
 \vdots & \cdots & \vdots \\
 \left[A_{n1}\right] & \ldots & \left[A_{nn}\right]\end{array}\right) \;. 
\end{equation}

In (\ref{e10})~ $C_i^T=c_i(x_1),c_i(x_2),\ldots ,c_i(x_m)$,~
 $G_i^T=g_i(x_1),g_i(x_2),\ldots ,g_i(x_m)$ and each block
 $[A_{ij}]$ is a tridiagonal matrix
$$
[A_{ij}]\! = \!\left(\!\! \begin{array}{ccccccc}
  \delta_{ij}\!+\!\frac{2}{h^2}a_{ij} \!&\! \frac{-1}{h^2}a_{ij} 
  \!&\! 0 \!&\! 0\!&\! \ldots \!&\! 0\!&\! 0 \\
  \frac{-1}{h^2}a_{ij} \!&\! \delta_{ij}\!+\!\frac{2}{h^2}a_{ij} 
  \!&\! \frac{-1}{h^2}a_{ij} \!&\! 0 \!&\! \ldots \!&\! 0 \!&\! 0 \\
  0 \!&\! \frac{-1}{h^2}a_{ij} \!&
  \! \delta_{ij}\!+\!\frac{2}{h^2}a_{ij} 
  \!&\! \frac{-1}{h^2}a_{ij} \!&\! \ldots \!&\! 0 \!&\! 0 \\
  \vdots \!&\!\vdots\!&\!\vdots\!&\!\vdots\!&\!\vdots\!&
  \!\vdots\!&\! \vdots \\
  0 \!&\! 0 \!&\! 0 \!&\! \ldots \!&\! \frac{-1}{h^2}a_{ij} \!&\! 
  \delta_{ij}\!+\!\frac{2}{h^2}a_{ij} \!&\! \frac{-1}{h^2}a_{ij}  \\
  0 \!&\! 0 \!&\! 0 \!&\! 0 \!&\! \ldots \!&
  \! \frac{-1}{h^2}a_{ij} \!&\! 
  \delta_{ij}\!+\!\frac{2}{h^2}a_{ij} 
\end{array}\! \right)\!. 
$$
In general ~$\underline{\mathcal{A}}$~ is real, non-symmetric
matrix (because ~$a_{ij}\neq a_{ji}$, see (\ref{e8})).

\subsection{Examples of numerical results for $d=1$.} 

Because solutions to equation (\ref{e1}) are traveling wave-like functions
we use free boundary conditions and large enough grid (precisely, 
~$f(x,t)_{x\rightarrow\pm\infty}{\longrightarrow}0$~
for any finite ~$t$). 

As initial condition we take a Gaussian distribution 
$ f(x,0) = g(x) = \exp \left[-\frac{x^2}{\sigma^2}\right]$.
It can represent the initial distribution of the temperature for 
 the heat equation ($\alpha =1$) or initial displacement of the
 medium for the wave equation  ($\alpha =2$).
 
For one spatial dimension $N=n\cdot m <10^4$ is usually sufficient for
obtaining a reasonable approximate numerical solution. For such
$N$ the set of linear equations (\ref{e9}) can be solved by standard
methods (e.g.\ LU decomposition). In fig.\ \ref{f1} we show 
numerical solutions to (\ref{e1}) for $\alpha=1, \frac{5}{4},
\frac{3}{2}, \frac{7}{4}$ and 2, at two particular time instants
$t=6$ and $t=12$. The value of $\sigma$ in the initial condition
was taken as $\sigma=1$.
The reader can easily see a transition from a diffusion-like solution
for $\alpha=1$, through intermediate cases for 
$1<\alpha <2$, to a wave-like solution for $\alpha=2$.

\begin{figure}[htb]  \label{f2}
\centerline{\resizebox{0.95\textwidth}{!}{\includegraphics{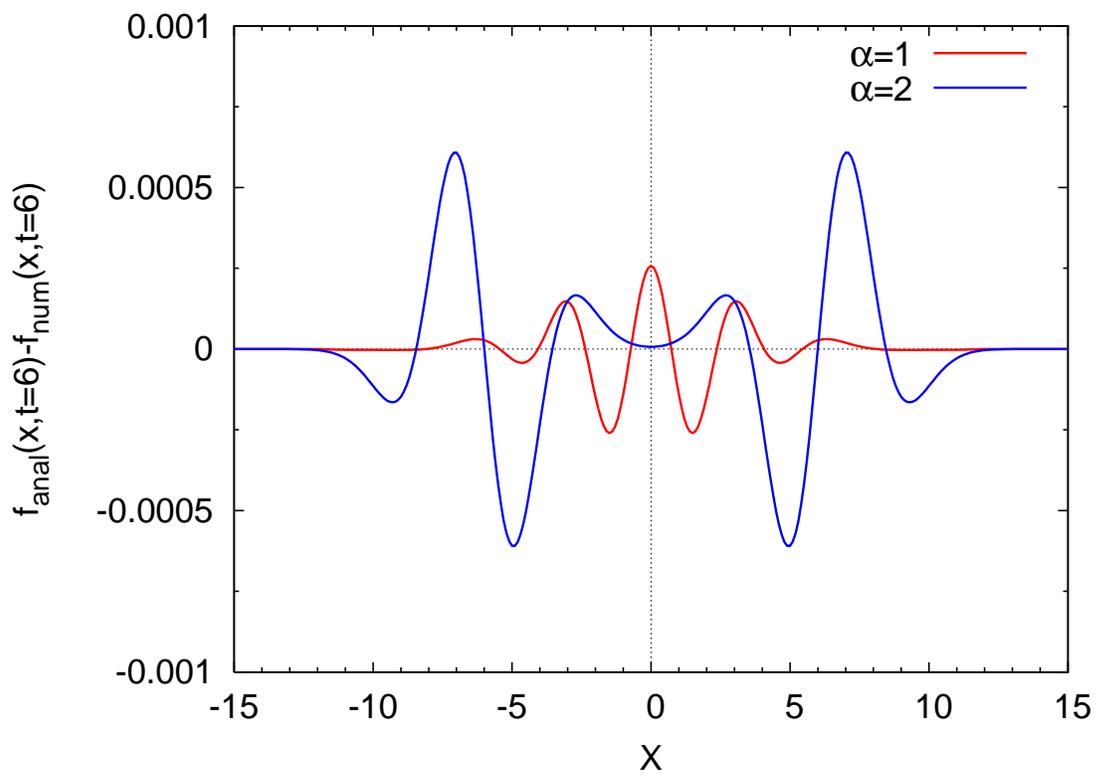}  }}
\caption{  
Numerical errors for $t=6$ and cases $\alpha=1$ (red line)
and $\alpha=2$ (blue line).	}    
\end{figure}

The knowledge of the analytical form of solutions for $\alpha=1$ and 2
allows us to keep approximation errors within a required range. 
To maintain the errors $\epsilon(x,t)=|f_\mathrm{anal}(x,t)-
f_\mathrm{num}(x,t)|\leq 10^{-3}$ it was enough, for $t=6$, to take
into account a grid of $m=151$ points in $x$-coordinate, covering
the interval $x\in [-15,15]$ and subspace $H_n$ with $n=8$.
For $t=12$ case and the same error bounds the grid had to be increased to 
$m=201$ points for the interval $x\in [-20,20]$
and subspace  $H_n$ to $n=18$.
Fig.\ \ref{f2} presents approximation errors $\epsilon(x,t)$ for
$t=6$ and $\alpha=1$ and 2.

For larger times the number of grid points and size of subspace $H_n$ has
to grow in order to keep the same precision of numerical solutions.
As the matrix $\underline{\mathcal{A}}$ is sparse
(among $n^2\cdot m^2$ elements of $\underline{\mathcal{A}}$ at most
$n^2 (3m-2)$ are non-zero) iterative methods for solving (\ref{e9})
become necessary.

\begin{figure}[hbt] 
\begin{center}
\resizebox{0.85\textwidth}{!}{\includegraphics{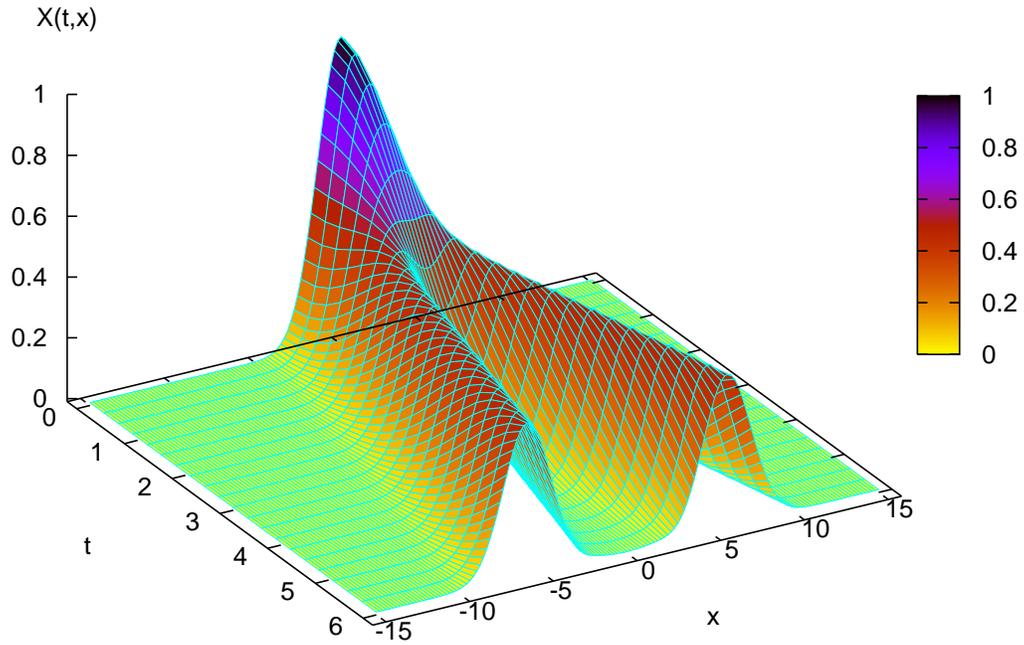}}
\resizebox{0.85\textwidth}{!}{\includegraphics{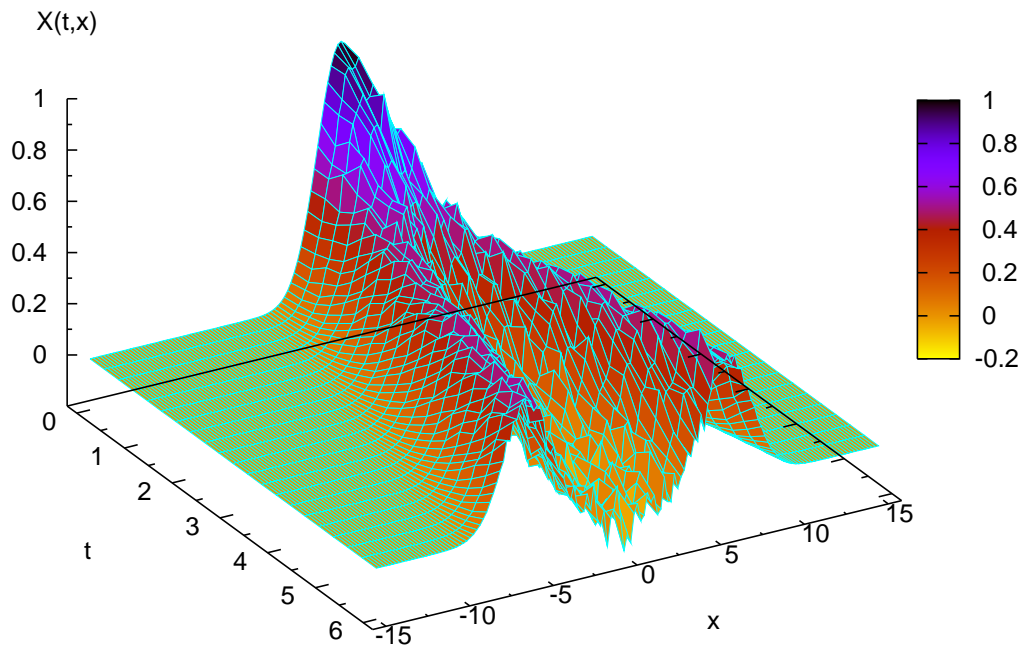}}
\end{center}
\caption{Time evolution of the solution to the equation (1) with 
$\alpha=2$ for $t\in [0,6]$, 
the deterministic solution (top) and the  solution to the equation (2),
a particular stochastic trajectory (bottom).}
\label{f2a}
\end{figure}

For stochastic equation (\ref{e4}) we need some assumptions for the 
process $W$. For the first attempt we assumed that the process $W$
is uniform in time, i.e.\ $W(x,t) =CW_1(x,t)$ (the constant $C$
represents a 'strength' of the stochastic forces). 
Then we can approximate
the convolution $\int_0^t S(t-s)dW(s,x)$ in (\ref{e4}) 
in the following way:
\begin{equation} \label{e15}
 \int_0^t S(t-s)dW(s,x)= 
 \sum_{i=0}^{I-1} S(t-s_i)[W(s_{i+1},x)-W(s_{i},x)],
\end{equation}
where the time interval $[0,t]$ was divided into a time grid 
$\{ t_i=i\tau,~i=0,1,\ldots,I\}$, $\tau=\frac{t}{I}$.
For cases $\alpha=1$ and 2, when $S(t)$ is known analytically
(see (\ref{e5a}) and (\ref{e5b})) 
the stochastic convolution can be computed
numerically. Fig.\ \ref{f2a} compares the time evolution of solutions 
obtained numerically for $\alpha=2$ and $t\in [0,6]$. 
The top part represents the solution of the deterministic equation
(\ref{e1}), the bottom one an example of a single stochastic trajectory
(solution of the stochastic equation  (\ref{e4}) with $W$ uniform 
in time, $C=0.1$).
For more details and examples of numerical results, see \cite{KaRo}.
 
\section{Galerkin method, case $d=2$.}

For  $d=2$ the equation (\ref{e9}) reads 
\begin{eqnarray} \label{e13}
 g_j(x_i,y_l) & = &  c_j(x_i,y_l) 
  + \frac{1}{h^2} \sum_{k=1}^n a_{jk}\, 
 [-c_k(x_{i-1},y_{l})-c_k(x_{i},y_{l-1}) \\ &&
   + 4 c_k(x_{i},y_{l}) 
 - c_k(x_{i+1},y_{l})- c_k(x_{i},y_{l+1})] \;, 
\end{eqnarray} 
with  $j=1,2,\ldots,n,~i=1,2,\ldots,m,~l=1,2,\ldots,m$.

Now, in matrix equation  
\begin{equation} \label{e13a}
\underline{\mathcal{A}}\, 
\underline{c} = \underline{g},   
\end{equation}
$\underline{c}~ $ and $~\underline{g}~$ are $N=n\cdot m^2$-dimensional
vectors, such that \hfill \\
$\underline{c}^T=c_1(x_1,y_1),\ldots,c_1(x_1,y_m), 
c_1(x_2,y_1) \ldots,c_1(x_2,y_m),\ldots,c_n(x_m,y_m)$, \hfill\\
$\underline{g}^T=g_1(x_1,y_1),\ldots,g_1(x_1,y_m), 
g_1(x_2,y_1) \ldots,g_1(x_2,y_m),\ldots,g_n(x_m,y_m)$
\hfill \\ and 
\begin{equation} \label{e13b}
\underline{\mathcal{A}} =  \left( \begin{array}{ccc}
 \left[A_{11}\right] & \ldots & \left[A_{1n}\right] \\ 
 \left[A_{21}\right] & \ldots & \left[A_{2n}\right] \\ 
 \vdots & \cdots & \vdots \\
 \left[A_{n1}\right] & \ldots & \left[A_{nn}\right]\end{array}\right).
\end{equation}
Now, every block $[A_{ij}]$ 
is the tridiagonal matrix composed of smaller blocks 
\begin{equation} \label{e14}
[A_{ij}] = \left( \begin{array}{ccccccc}
 (\alpha)_{ij}  & (\beta)_{ij} &  (0)  &  (0) &  (0) & \cdots & (0) \\
 (\beta)_{ij}  & (\alpha)_{ij} &  (\beta)_{ij}  &  (0) &  (0) & \cdots & (0) \\ 
 (0)  & (\beta)_{ij}  & (\alpha)_{ij}  &  (\beta)_{ij} &  (0) & \cdots & (0) \\
 \vdots  &  \cdots &  &  \cdots  &   & \cdots&   \vdots \\
 (0) & (0) & \cdots & (0) & (\alpha)_{ij} & (\beta)_{ij} & (0) \\
 (0) & (0) & (0) & \cdots & (\beta)_{ij} & (\alpha)_{ij} & (\beta)_{ij} \\
 (0) & (0) & (0) & (0) & \cdots & (\beta)_{ij} & (\alpha)_{ij} 
\end{array} \right)  .
\end{equation}
Blocks $(\alpha_{ij})$ are  tridiagonal 
\begin{equation} \label{e14a}
 \!\!(\alpha_{ij})\!\! = \!\!\left(\!\!\! \begin{array}{ccccccc}
  \delta_{ij}\!+\!\frac{4}{h^2}a_{ij} \!\!&\!\! \frac{-1}{h^2}a_{ij} 
  \!\!&\!\! 0 \!\!&\!\! 0\!\!&\!\! \ldots \!\!&\!\! 0\!\!&\!\! 0 \\
  \frac{-1}{h^2}a_{ij} \!\!&\!\! \delta_{ij}\!+\!\frac{4}{h^2}a_{ij} 
  \!\!&\!\! \frac{-1}{h^2}a_{ij} 
  \!\!&\!\! 0 \!\!&\!\! \ldots \!\!&\!\! 0 \!\!&\!\! 0 \\
  0 \!\!&\!\! \frac{-1}{h^2}a_{ij} \!\!&\!\! \delta_{ij}\!+\!\frac{4}{h^2}a_{ij} 
  \!&\! \frac{-1}{h^2}a_{ij} \!\!&\!\! \ldots \!\!&\!\! 0 \!\!&\!\! 0 \\
  \vdots \!\!&\!\!\vdots\!\!&\!\!\vdots\!\!&\!\!\vdots\!\!&
  \!\!\vdots\!\!&\!\!\vdots\!\!&\!\! \vdots \\
  0 \!\!&\!\! 0 \!\!&\!\! 0 \!\!&\!\! \ldots \!\!&\!\! \frac{-1}{h^2}a_{ij} \!&\! 
  \delta_{ij}\!+\!\frac{4}{h^2}a_{ij} \!\!&\!\! \frac{-1}{h^2}a_{ij}  \\
  0 \!\!&\!\! 0 \!\!&\!\! 0 \!\!&\!\! 0 \!\!&\!\! \ldots 
  \!\!&\!\! \frac{-1}{h^2}a_{ij} \!\!&\!\! 
  \delta_{ij}\!+\!\frac{4}{h^2}a_{ij} 
\end{array}\!\! \right) \!\! , 
\end{equation}
blocks $(\beta_{ij})$ are  diagonal 
\begin{equation}  \label{e16}
(\beta_{ij}) = \left( \begin{array}{ccccccc}
  \frac{-1}{h^2}a_{ij} & 0 & 0 & 0 & 0 & \ldots & 0 \\
  0 & \frac{-1}{h^2}a_{ij} & 0 & 0 & 0 & \ldots & 0 \\
  0 & 0 & \frac{-1}{h^2}a_{ij} & 0 & 0 & \ldots & 0 \\
  \vdots & \vdots & \vdots & \vdots & \vdots & \vdots & \vdots \\
  0 & 0 & 0 & \ldots & 0 & \frac{-1}{h^2}a_{ij} & 0 \\
  0 & 0 & 0 & 0 & \ldots & 0 & \frac{-1}{h^2}a_{ij} 
\end{array}\right)
\end{equation}
and (0) are zeros, each of size  $m\cdot m$.

Dimension of vectors ~$\underline{c}~ \mbox{and} ~\underline{g}$
is ~$N=n\cdot m^2$. Already for 
$n\!=\!16,~m\!=\!200$, $N$ becomes large, reaching value of  
$N=6.4\cdot 10^5$
unknowns and the number of matrix elements for the matrix  
$\underline{\mathcal{A}}$ reaches $~N^2=n^2m^4 >4\cdot 10^{11}$.

Fortunately, the matrix $\underline{\mathcal{A}}$ is sparse. Blocks 
~$[A_{ij}]$ have at most ~$m(3m-2)+2(m-1)m$ non-zero elements.
Then the number of non-zero elements of matrix 
$\underline{\mathcal{A}}$ is at most (some ~$a_{ij}$ could be
$0$)
\begin{equation}  \label{e18}
\mathcal{N}\leq n^2m(5m-4).
\end{equation}
For ~$n=16,~m=200$, 
\begin{equation}\mathcal{N}\leq 5.12\cdot 10^{7}
\end{equation}

The size of ~$\underline{\mathcal{A}}$ and its
sparseness property makes using {\it iterative methods} 
for $d\geq 2$ necessary.

\subsection{Conjugate and Bi-Conjugate Gradient Method }

If matrix $A$ is symmetric and positive definite, then the problem 
$A\underline{x}=\underline{b}$ is equivalent to minimizing the function
$f(\underline{x})=\frac{1}{2}\underline{x}\cdot A\cdot
\underline{x}-\underline{b}\cdot\underline{x}$.
This function is minimized when its gradient $\nabla f =A\cdot
\underline{x}-\underline{b}$ is zero. 
In such case a variant of  {\it Conjugate Gradient Method}~ is
applicable \cite{NumR,Temp,Shew}.
 
In our problems the matrix $\underline{\mathcal{A}}$ is {\it not} symmetric. 
Therefore a more complicated {\it Bi-Conjugate Gradient Method}~ 
is required \cite{NumR,Temp,Shew}.

\subsection{Preconditionig } 

The convergent rate of iterative methods depends strongly on spectral 
properties of the matrix ~$\underline{\mathcal{A}}$.
Usually matrix ~$\underline{\mathcal{A}}$ is 
{\it ill-conditiond}. The condition number 
$\kappa = \lambda_{\mathrm{max}}/ \lambda_{\mathrm{min}}$ is big ($\lambda_{i}$
denotes an  eigenvalue). Then the convergence of iterations is usually
so slow that accumulation of numerical errors often makes obtaining the
solution impossible.
The remedy is {\it preconditioning}. Suppose that $M$ is a matrix
that approximates  ~$\underline{\mathcal{A}}$, but easier to invert.
We can solve $\underline{\mathcal{A}}\cdot \underline{x}=\underline{b}$~
indirectly by solving
~$M^{-1}\cdot\underline{\mathcal{A}}\cdot \underline{x} =
M^{-1}\cdot\underline{b}$.
If $\kappa(M^{-1}\underline{\mathcal{A}})\ll\kappa(\underline{\mathcal{A}})$,
the number of iterations is reduced significantly.

There are several ways of choosing a preconditioner matrix $M$.
In our case we can take an advantage of knowing detailed structure
of matrix $\underline{\mathcal{A}}\,$ (\ref{e13b}-(\ref{e14a})
which all elements are related
to the elements of small matrix $\underline{a}$ (\ref{e8}).
Blocks $[A_{ij}]$ (\ref{e14}) and $(\alpha_{ij})$ (\ref{e14a})
are tridiagonal. We choose the preconditioner matrix $M$
in the same block form as the matrix $\underline{\mathcal{A}}\,$,
but leaving only diagonal blocks $(\alpha_{ij})$ in (\ref{e14})
and diagonal elements $\delta_{ij}+\frac{1}{h^2}a_{ij}$
in (\ref{e14a}). All other elements of $(\alpha_{ij})$  and 
$(\beta_{ij})$ are set equal zero.
Then the matrix $M$ has block form with diagonal blocks containing
the same element $\gamma_{ij}=\delta_{ij}+\frac{4}{h^2}a_{ij}$
on their diagonals.
Hence the matrix $M^{-1}$ has the same block structure,
with elements $\gamma^{-1}$ on block's diagonals. 
The size of $\gamma$ is only $n\cdot n$, so $\gamma^{-1}$
can be calculated easily by standard methods with the machine precision.
The resulting matrix $M^{-1}\underline{\mathcal{A}}\,$
has usually the condition number several orders of magnitude smaller
than that of the original matrix $\underline{\mathcal{A}}\,$.
In calculations leading to results presented below 
such kind of preconditionig
allows to obtain a reasonable accuracy within $10^2-10^4$ iterations
for problems with $\sim 10^5$ unknowns. 

\subsection{Numerical results}

\begin{figure}[htb]  
\centerline{\resizebox{0.85\textwidth}{!}{\includegraphics{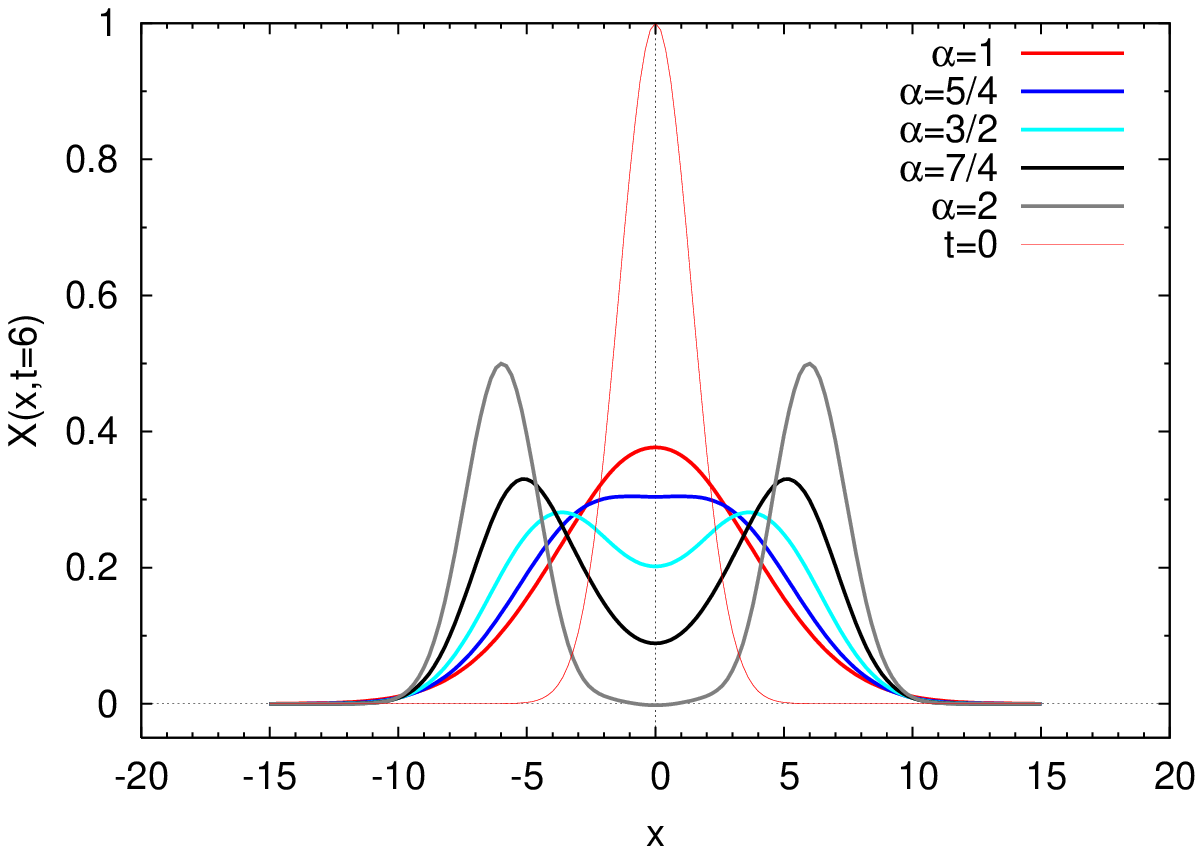}} }
\centerline{\resizebox{0.85\textwidth}{!}{\includegraphics{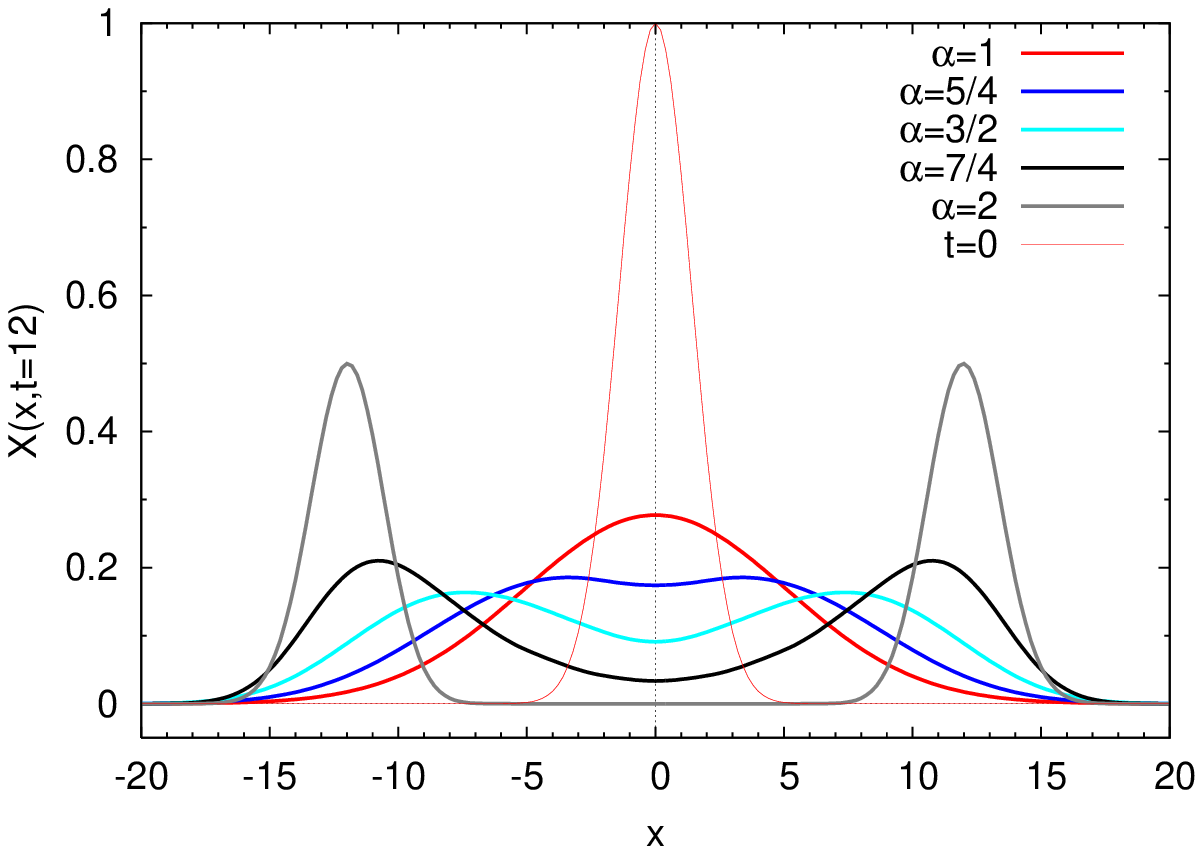}  } }
\caption{ 
Numerical solutions to deterministic Volterra
equation (\ref{e1}) at $t=6$ (top) and $t=12$ (bottom). Presented 
are cases with $\alpha\!=\!1,\frac{5}{4},\frac{3}{2},\frac{7}{4},2$. 
Thin red line represents $f(x,t=0)$ with $\sigma=2$.  	}  
\label{f1} 
\end{figure}

\begin{figure}[htb]
\centerline{\resizebox{0.95\textwidth}{!}{\includegraphics{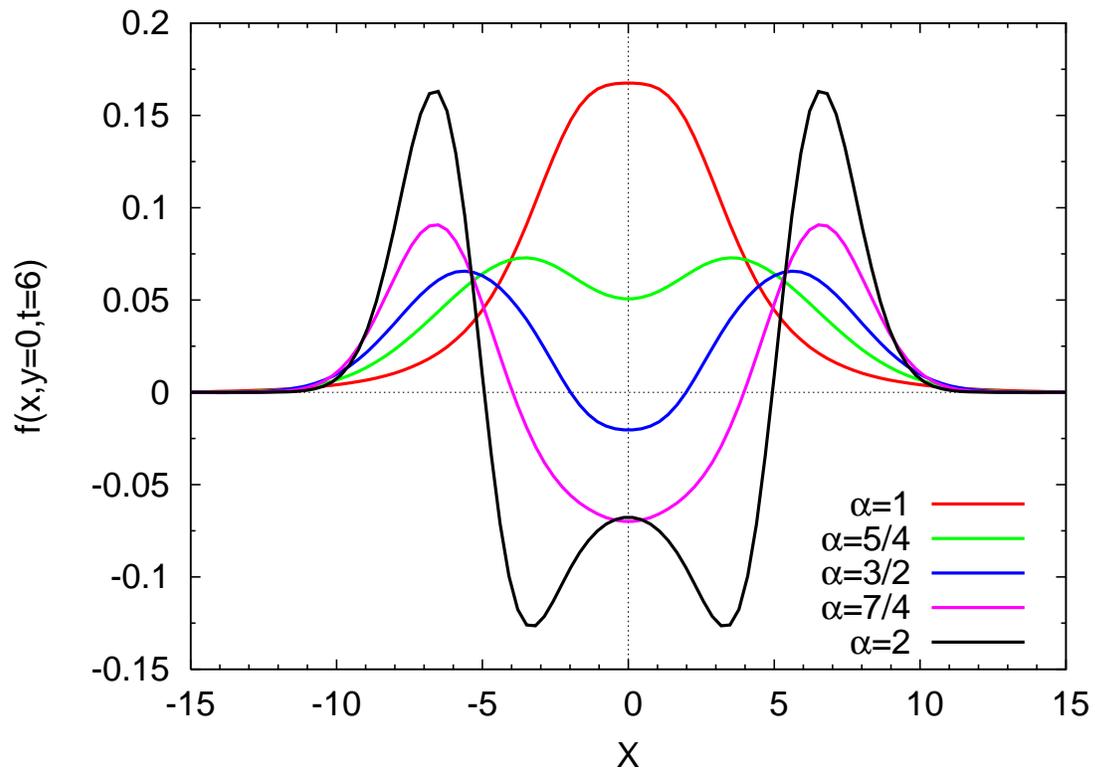}  }}
\caption{ \small 
One dimensional sections of two dimensional 
numerical solutions to deterministic Volterra
equation (\ref{e1}) at ~$t=6$. Presented 
are cases with ~$\alpha=1,\frac{5}{4},\frac{3}{2},\frac{7}{4},2$
for radially symmetric initial conditions.  }  
\label{f3}  
\end{figure}

Solutions to deterministic equation (\ref{e1}) for $d=1$ and $d=2$
differ substantially from each other. In fig.\ \ref{f3} we present
the numerical solutions to (\ref{e1}) for 2 spatial dimensions
in a way convenient for comparison with fig.\  \ref{f1} (top)
presenting solutions to (\ref{e1}) for 1 spatial dimension.
The initial condition for results displayed in fig.\  \ref{f3} 
is in the form $f(x,y,0)=\exp [-\frac{x^2+y^2}{\sigma^2}],~\sigma=2$. 
The curves in  fig.\  \ref{f3} represent cuts of solutions
$f(x,y,t)$ along $y=0$, i.e.~$f(x,0,t)$ and it is clearly seen that for
all given values of $\alpha$ the profiles of the solutions for
$d=1$ and $d=2$ are different.

\begin{figure}[htb] \
\centerline{
\resizebox{0.95\textwidth}{!}{\includegraphics{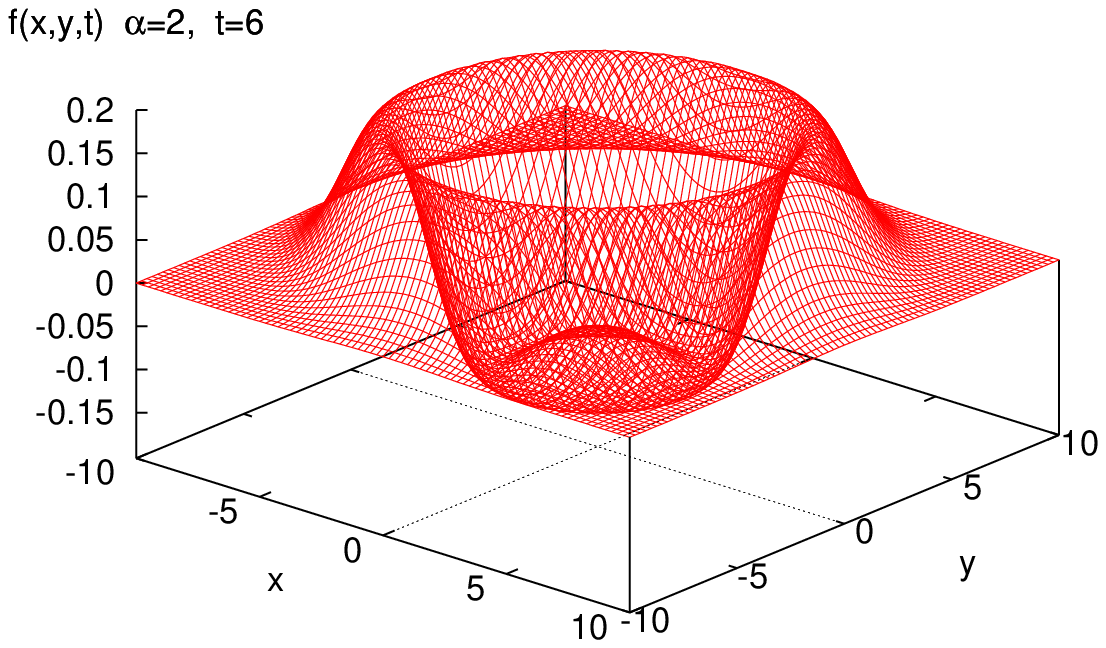}  } }
\centerline{
\resizebox{0.95\textwidth}{!}{\includegraphics{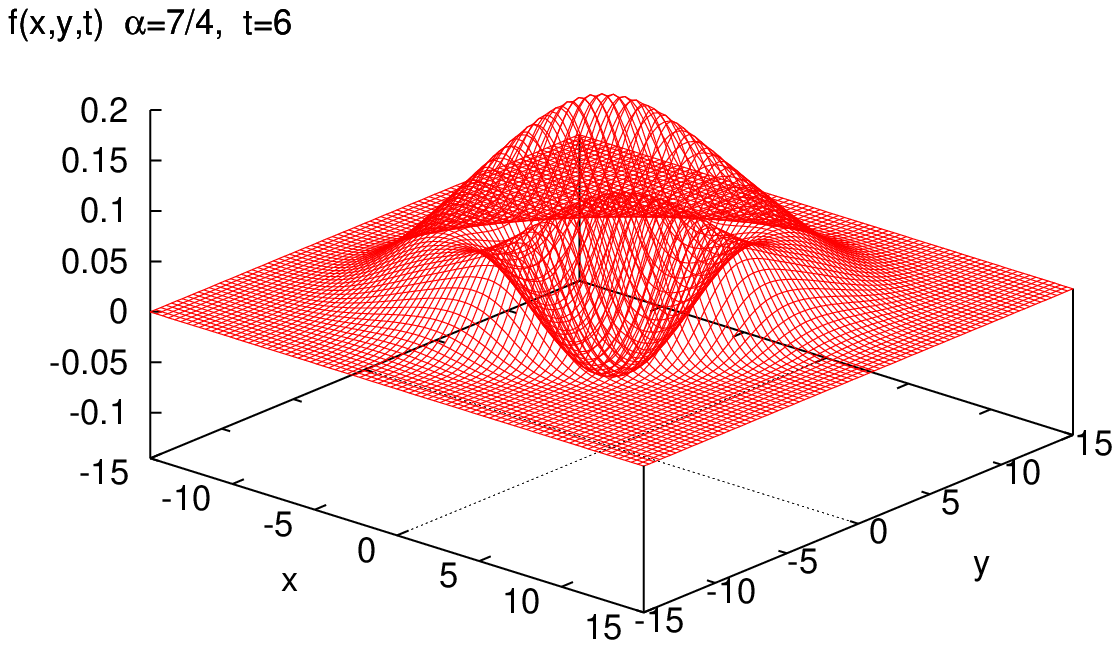}  }}
\caption{ \small 
Examples of numerical solutions to deterministic Volterra
 equation (\ref{e1}) for2 spatial dimensions. Presented 
are cases with $\alpha=2$, radially symmetric initial 
condition (top) at $t=6$, and $\alpha=\frac{7}{4}$, asymmetric initial
condition (bottom) at $t=6$.}    
 \label{f4}
\end{figure}

In fig.\ \ref{f4} two examples of the solutions $f(x,y,t)$ at $t=6$ are
displayed. In the upper part the case $\alpha=2$ and radially symmetric
initial condition is shown.
In the lower part the case   $\alpha=\frac{7}{4}$ with radially
asymmetric  initial condition 
($f(x,y,0)=\exp [-\frac{(x+y)^2}{\sigma_1^2}
-\frac{(x-y)^2}{\sigma_2^2}]$, with $\sigma_1=4$, $\sigma_2=2$) 
is presented.

\end{document}